\documentclass[12pt]{article}
\usepackage{graphicx}
\usepackage{amsmath,amsthm,amssymb,enumerate}
\usepackage{euscript,mathrsfs}
\usepackage{color}
\usepackage{dsfont}
\usepackage[notref,notcite]{showkeys}
\usepackage[left=2cm,right=2cm,top=3.5cm,bottom=3.5cm]{geometry}
\usepackage{color}
\usepackage[framemethod=tikz]{mdframed}
\usepackage{bm}
\allowdisplaybreaks

\usepackage{esint}
\usepackage{soul}

\catcode`\@=11 \@addtoreset{equation}{section}

\catcode`\@=12

\allowdisplaybreaks

\newtheorem{Theorem}{Theorem}[section]
\newtheorem{Proposition}[Theorem]{Proposition}
\newtheorem{Lemma}[Theorem]{Lemma}
\newtheorem{Corollary}[Theorem]{Corollary}

\theoremstyle{definition}
\newtheorem{Definition}[Theorem]{Definition}

\newtheorem{Remark}[Theorem]{Remark}

\newcommand{\bTheorem}[1]{
\begin{Theorem} \label{T#1} }
\newcommand{\eT}{\end{Theorem}}

\newcommand{\bProposition}[1]{
\begin{Proposition} \label{P#1}}
\newcommand{\eP}{\end{Proposition}}

\newcommand{\bLemma}[1]{
\begin{Lemma} \label{L#1} }
\newcommand{\eL}{\end{Lemma}}

\newcommand{\bCorollary}[1]{
\begin{Corollary} \label{C#1} }
\newcommand{\eC}{\end{Corollary}}

\newcommand{\bRemark}[1]{
\begin{Remark} \label{R#1} }
\newcommand{\eR}{\end{Remark}}

\newcommand{\bDefinition}[1]{
\begin{Definition} \label{D#1} }
\newcommand{\eD}{\end{Definition}}

\newcommand{\Del}{\Delta_x}

\newcommand{\Ds}{\mathbb{D}_x}

\newcommand{\bfphi}{\boldsymbol{\varphi}}

\newcommand{\bFormula}[1]{
\begin{equation} \label{#1}}
\newcommand{\eF}{\end{equation}}

\newcommand{\Ov}[1]{\overline{#1}}

\newcommand{\aleq}{\stackrel{<}{\sim}}

\newcommand{\vr}{\varrho}
\newcommand{\vre}{\vr_\ep}

\newcommand{\vue}{\vu_\ep}

\newcommand{\vu}{\vc{u}}

\newcommand{\vc}[1]{{\bm #1}}

\newcommand{\Div}{{\rm div}_x}
\newcommand{\Grad}{\nabla_x}

\newcommand{\dx}{\,{\rm d} {x}}
\newcommand{\dz}{{\rm d} {z}}
\newcommand{\dt}{\,{\rm d} t }

\newcommand{\intO}[1]{\int_{\Omega} #1 \ \dx}
\newcommand{\avintBe}[1]{\fint_{B_{\ep, t}} #1 \dx}

\newcommand{\intRd}[1]{\int_{R^3} #1 \ \dx}

\newcommand{\D}{{\rm d}}

\newcommand{\ep}{\varepsilon}

\newcommand{\br}{ \nonumber \\ }

\def\softd{{\leavevmode\setbox1=\hbox{d}%
          \hbox to 1.05\wd1{d\kern-0.4ex{\char039}\hss}}}
\definecolor{Cgrey}{rgb}{0.85,0.85,0.85}
\definecolor{Cblue}{rgb}{0.50,0.85,0.85}
\definecolor{Cred}{rgb}{1,0,0}
\definecolor{fancy}{rgb}{0.10,0.85,0.10}

\newcommand\Cbox[2]{%
    \newbox\contentbox%
    \newbox\bkgdbox%
    \setbox\contentbox\hbox to \hsize{%
        \vtop{
            \kern\columnsep
            \hbox to \hsize{%
                \kern\columnsep%
                \advance\hsize by -2\columnsep%
                \setlength{\textwidth}{\hsize}%
                \vbox{
                    \parskip=\baselineskip
                    \parindent=0bp
                    #2
                }%
                \kern\columnsep%
            }%
            \kern\columnsep%
        }%
    }%
    \setbox\bkgdbox\vbox{
        \color{#1}
        \hrule width  \wd\contentbox %
               height \ht\contentbox %
               depth  \dp\contentbox
        \color{black}
    }%
    \wd\bkgdbox=0bp%
    \vbox{\hbox to \hsize{\box\bkgdbox\box\contentbox}}%
    \vskip\baselineskip%
}

\mdfdefinestyle{MyFrame}{%
	linecolor=black,
	outerlinewidth=1pt,
	roundcorner=5pt,
	innertopmargin=\baselineskip,
	innerbottommargin=\baselineskip,
	innerrightmargin=10pt,
	innerleftmargin=10pt,
	backgroundcolor=white!20!white}

\date{}

\makeindex
\begin{document}


\title{On the motion of a small rigid body in a viscous compressible fluid}

\author{Eduard Feireisl
	\thanks{The work of E.F. was partially supported by the
		Czech Sciences Foundation (GA\v CR), Grant Agreement
		21--02411S. The Institute of Mathematics of the Academy of Sciences of
		the Czech Republic is supported by RVO:67985840. A.R and A.Z have been partially supported by the Basque Government through the BERC 2022-2025 program and by the Spanish State Research Agency through BCAM Severo Ochoa excellence accreditation SEV-2017-0718 and through project PID2020-114189RB-I00 funded by Agencia Estatal de Investigación (PID2020-114189RB-I00 / AEI / 10.13039/501100011033). A.Z. was also partially supported  by a grant of the Ministry of Research, Innovation and Digitization, CNCS - UEFISCDI, project number PN-III-P4-PCE-2021-0921, within PNCDI III.} 
	\and Arnab Roy$^1$ \and Arghir Zarnescu$^{1,2,3}$
}

\date{\today}

\maketitle

\bigskip

\centerline{$^*$  Institute of Mathematics of the Academy of Sciences of the Czech Republic}

\centerline{\v Zitn\' a 25, CZ-115 67 Praha 1, Czech Republic}

\centerline{$^1$ BCAM, Basque Center for Applied Mathematics}

\centerline{Mazarredo 14, E48009 Bilbao, Bizkaia, Spain}

\centerline{$^2$IKERBASQUE, Basque Foundation for Science, }

\centerline{Plaza Euskadi 5, 48009 Bilbao, Bizkaia, Spain}

\centerline{$^3$`Simion Stoilow" Institute of the Romanian Academy,}

\centerline{21 Calea Grivi\c{t}ei, 010702 Bucharest, Romania }

\maketitle

\begin{abstract}
	
	We consider the motion of a small rigid object immersed in a viscous compressible fluid in the 3-dimensional Eucleidean space. Assuming the object is a ball of a small radius $\ep$ we show that the behavior of the 
	fluid is not influenced by the object in the asymptotic limit $\ep \to 0$. The result holds for the isentropic pressure law $p(\vr) = a \vr^\gamma$ for any $\gamma > \frac{3}{2}$ under 
	mild assumptions concerning the rigid body density. In particular, the latter may be bounded as soon as $\gamma > 3$.  
	The proof uses a new method of construction of the test functions in the weak formulation of the problem, and, in particular, a new form of the so-called Bogovskii operator.

\end{abstract}

\bigskip


\section{Introduction}
\label{i}

Consider a rigid body immersed in a viscous fluid. Intuitively, the impact of a ``small'' body on the fluid motion should be negligible. A rigorous justification of this statement has 
been obtained in several recent studies on condition that the fluid is incompressible, see 
Lacave and Takahashi \cite{LacTak}, Iftimie et al. \cite{MR2244381}, He and Iftimie \cite{HeIft1, HeIft2}, Dashti and Robinson \cite{MR2781594}. The approach of 
\cite{LacTak} is based on the $L^p-L^q$ estimates for the associated solution semigroup available in the 
2d-setting, while He and Iftimie \cite{HeIft1} use a specific construction of time dependent test functions vanishing on the moving body. In \cite{MR2781594}, a viscous fluid-rigid disc system has been studied where the disc is not rotating and they proved that the body does not influence the flow
in the asymptotic limit. Lacave \cite{MR2557320} studies the limit of a viscous fluid flow
in the exterior of a thin obstacle shrinking to a curve. In \cite{FRZ22}, the authors consider the motion of a rigid body inside a compressible fluid in planar domain and establish that the influence of the body
on the fluid is negligible if the diameter of the body is small and the fluid is nearly incompressible (the low Mach number regime).

 Recently, Bravin and Ne\v casov\' a \cite{BraNec} combined the 
technique of \cite{HeIft1} with the pressure estimates obtained via the new Bogovskii operator introduced in \cite{DieFeLu} and Lu and Schwarzacher \cite{LuSchw} to handle the 3d compressible case under certain technical restrictions imposed on the pressure--density equation of state, notably on the value of adiabatic exponent. The above mentioned technique seems difficult to adapt to the planar (2d) motion of a compressible fluid and the results are not optimal even in the 3d-setting, where certain additional 
restrictions are needed on the value of the adiabatic exponent.  
Indeed a single point in the $d$-dimensional space has a positive $W^{1,p}-$capacity as soon as $p > d$. Accordingly, 
the approximation technique developed in \cite{HeIft1} requires the pressure to be uniformly 
$\frac{d}{d-1}$ integrable when the diameter of the body approaches zero. Unfortunately, the best known estimates for the standard example of the isentropic pressure $p(\vr) = a \vr^{\gamma}$ read 
\[
p(\vr) \in L^q,\ \mbox{with}\ q = \frac{d + 2}{d} - \frac{1}{\gamma} 
\]  
see Lions \cite{LI4}, meaning the value $q = 2$ for $d=2$ is never achieved, while $q = \frac{3}{2}$ for $d = 3$ requires $\gamma \geq 6$.

To handle physically realistic adiabatic exponents, we propose a new approach based on the concept of weak solution introduced in \cite{EF64}. We first observe that the test functions used for the approximate problem need not vanish on the moving body but only 
satisfy the rigid body motion constrain. Using this rather straightforward observation we construct a new 
approximation operator based on the version of the Bogovskii operator on uniformly John domains due to 
Diening, R\accent23u\v zi\v cka, and Schumacher \cite{DieRuzSch}. The result seems optimal as we recover the desired convergence without any additional restrictions on the equation of state, notably on the adiabatic 
coefficient $\gamma > \frac{3}{2}$, $d=3$ in agreement with the available existence theory.

The paper is organized as follows. In Section \ref{P}, we formulate the problem, recall the concept of weak solution and state the main result of the paper. The available uniform bounds are summarized in Section \ref{U}. Sections \ref{T} and \ref{E} are the heart of the paper. We construct a general restriction operator along with its vector valued version preserving 
the divergence of the extended function. The pressure estimates necessary to perform the asymptotic limit 
for ``vanishing'' body are derived in Section \ref{pe}. Finally, the convergence proof is completed in Section 
\ref{c}.

\section{Problem formulation, weak solutions, main results}
\label{P}

The motion of a compressible viscous fluid in the barotropic regime is governed by 

\begin{mdframed}[style=MyFrame]
	
{\bf Navier--Stokes system} 	
	
	\begin{align}
	\partial_t \vr + \Div (\vr \vu) &= 0, \label{P1}\\
	\partial_t (\vr \vu) + \Div (\vr \vu \otimes \vu) + \Grad p(\vr) &= \Div \mathbb{S}(\Grad \vu), 
	\label{P2}
	\end{align}
supplemented with

\noindent
{\bf Newton's rheological law} 
\begin{equation} \label{P3}
	\mathbb{S}(\Grad \vu) = \mu \left( \Grad \vu + \Grad^t \vu - \frac{2}{d} \Div \vu \mathbb{I} \right) + 
	\eta \Div \vu \mathbb{I},\ \mu > 0,\ \eta \geq 0.
	\end{equation}

	\end{mdframed}

\noindent
Here, $\vr$ is the mass density and $\vu$ is the fluid velocity. 

For mostly technical reasons, we 
focus on the Cauchy problem for $d = 3$ and neglect the effect of external forces. Accordingly, the fluid occupies the whole physical space $R^3$, where the density and the velocity satisfy the 
far field conditions
\begin{equation} \label{P4} 
\vu \to 0,\ \vr \to 0 \ \mbox{as}\ |x| \to \infty.
\end{equation}
In particular, we suppose the total mass of the fluid--body system is finite, 
\[
\intRd{ \vr(t, \cdot) } < \infty.
\]
More general far field conditions 
\[
\vu \to 0,\ \vr \to \vr_\infty \ \mbox{as}\ |x| \to \infty,\ \vr_\infty \geq 0 - \mbox{constant,}
\]
can be handled in a similar fashion.

We suppose the rigid body is a ball of the radius $\ep$ occupying at a given time $t \geq 0$ the compact set 
\[
B_{\ep,t} = \left\{ x \in R^3 \ \Big|\ |x - \vc{h}_\ep(t) | \leq \ep \right\}.
\]
We suppose that the mass density of the body $\vr_{\ep, B} > 0$ is a positive constant and the motion of the body is determined by the rigid velocity field
\[
\vu_{\ep,B} (t,x) = \vc{Y}_\ep (t) + \mathbb{Q}_{\ep,t}(t) (x - \vc{h}_\ep(t)),\ \frac{\D }{\dt} \vc{h}_\ep (t) = \vc{Y}_{\ep}(t).
\]
Accordingly, the fluid domain $Q_f$ is defined as
\[
Q_f = [0, T) \times R^3 \setminus \cup_{t \in [0,T)} B_{\ep,t} \subset 
[0,T) \times R^3.
\]

\subsection{Weak solutions}

Following \cite{EF64} we introduce a concept of \emph{weak solution} of the fluid--body interaction problem.

\begin{mdframed}[style=MyFrame]

\begin{Definition} [{\bf Weak solution}] \label{PD1}
	
	We say that $(\vre, \vue)$ is \emph{weak solution} of the fluid--body interaction problem with the initial 
	state $\vre(0, \cdot) = \vr_{\ep,0}$, $\vre \vue (0, \cdot) = \vc{q}_{\ep, 0}$ if the following holds:
	
	\begin{itemize}
	\item {\bf Compatibility.} $\vre \in L^\infty(0,T; L^1 \cap L^\gamma (R^3))$, 
	\[
	\vre(t,x) = \left\{ \begin{array}{l} \vr_{\ep,B} \ \mbox{if}\ x \in B_{\ep,t}, \\  \geq 0 \ \mbox{otherwise}      \end{array} \right. ,
	\]
	\[
	\intRd{ \vre (t, \cdot) } = \intRd{ \vr_{\ep,0} } \ \mbox{for any}\ t \in [0,T);
	\]
	$\vue \in L^2(0,T; D^{1,2}(R^3; R^3))$, 
	\[
	\vue(t,x) =  \vu_{\ep, B}(t,x) \ \mbox{if}\ x \in  B_{\ep,t} ;
	\]

\item {\bf Equation of continuity.} 
The integral identity 
\begin{equation} \label{P5}
\int_0^T \intRd{ \Big[ \vre \partial_t \varphi + \vre \vue \cdot \Grad \varphi \Big] } \dt = - \intRd{ \vr_{0, \ep} \varphi (0, \cdot) }
\end{equation}
holds for any $\varphi \in C^1_c([0,T) \times R^3)$. In addition, the renormalized equation
\begin{align}
\int_0^T &\intRd{ \Big[ b(\vre) \partial_t \varphi + b(\vre) \vue \cdot \Grad \varphi + ( b(\vre) - b'(\vre) \vre ) \Div \vue \varphi \Big] } \dt \br &= - \intRd{ b(\vr_{0, \ep}) \varphi (0, \cdot) }
 \label{P6}
\end{align}
holds for any $\varphi \in C^1_c([0,T) \times R^3)$ and any $b \in C^1[0, \infty)$, $b' \in C_c[0, \infty)$.

\item {\bf Momentum equation.} 
The integral identity
\begin{align}
\int_0^T &\intRd{ \Big[ \vre \vue \cdot \partial_t \bfphi + \vre \vue \otimes \vue : \Grad \bfphi + p(\vre) \Div \bfphi  } \dt
\br &= \int_0^T \intRd{ \mathbb{S}(\Grad \vue) : \Grad \bfphi } \dt - \intRd{ \vc{q}_{0, \ep} \cdot \bfphi (0, \cdot) }
\label{P7}
\end{align}
holds for any $\bfphi \in C^1_c ([0,T) \times \Omega; R^3)$ such that 
\begin{equation} \label{P8}
\Ds \bfphi (t,\cdot ) \equiv 
\frac{1}{2} \left( \Grad \bfphi + \Grad^t \bfphi \right)(t,\cdot) = 0 \ \mbox{on an open neighborhood of}\ B_{\ep,t}.
\end{equation}

\item {\bf Energy inequality.}
\begin{align} 
&\intRd{ \frac{1}{2} \vre |\vue|^2 (\tau, \cdot) } + 
\int_{R^3 \setminus B_{\ep, \tau}} P(\vre)  (\tau, \cdot) \dx  + \int_0^\tau \intRd{ 
	\mathbb{S} (\Grad \vue) : \Grad \vue } \dt \br &\quad \leq \intRd{ \frac{1}{2} \frac{ |\vc{q}_{0, \ep}|^2 }{\vr_{0, \ep}   } } + 
\int_{R^3 \setminus B_{\ep, 0}} P(\vr_{0, \ep}) \dx
\label{P9}
\end{align}
for a.a. $\tau \in (0,T)$, where
\[
P'(\vr) \vr - P(\vr) = p(\vr) \mbox{ or equivalently }P(\vr)=\vr\int\limits_1^{\vr} \frac{p(\tau)}{\tau^2}\ d\tau.
\]

\end{itemize}

\end{Definition}

\end{mdframed}

\begin{Remark} \label{HSS}
	
	The homogeneous Sobolev space $D^{1,2}(R^3)$ is defined as 
	\[
	D^{1,2}(R^3) = \left\{ \vc{v} \in L^6 (R^3) \ \Big| \  \Grad \vc{v} \in L^2(R^3) \right\}.
	\]
	
\end{Remark}

The existence of global--in--time weak solutions under the hypothesis $p \approx \vr^\gamma$,
$\gamma > \frac{3}{2}$ in a \emph{bounded} domain $\Omega \subset R^3$ was proved in \cite[Theorem 4.1]{EF64}.
The extension to the present setting is straightforward. The form of 
the energy inequality\eqref{P9} follows from \cite[formula (2.6) and Lemma 3.2]{EF64}.

\subsection{Main result}

Let us denote 
\[
\vr_{\ep, f}(t, \cdot ) = \vr_\ep (t, \cdot ) \mathds{1}_{R^3 \setminus B_{\ep, t}} 
\]
the fluid density.
We are ready to state our main result.

\begin{mdframed}[style=MyFrame]
	
	\begin{Theorem} [{\bf Convergence}] \label{PT1}
		
		Let the pressure $p$ be given by the isentropic equation of state
		\[ 
		p(\vr) = a \vr^\gamma,\ a > 0,\ \gamma > \frac{3}{2}.
		\]
		Let the density of the rigid body $\vr_{\ep,B}$ be a positive constant satisfying 
		\begin{align} 
			\vr_{\ep,B} \geq \underline{\vr} > 0,\ \ep^{-\underline{\beta}} \aleq \vr_{\ep,B} \aleq \ep^{- \Ov{\beta}} \ \mbox{as}\ \ep \to 0\br  \mbox{for some}\ 2 \left(\frac{3 - \gamma}{\gamma}\right) < \underline{\beta} 
			\leq \Ov{\beta} < 2.
			\label{P10}
			\end{align}	
		Finally, suppose that the initial data and energy satisfy
		\begin{align} 
	\vr_{0, \ep} &> 0,\ \vr_{0, \ep} \to \vr_0 \ \mbox{weakly in}\ L^1(R^3),\  \vc{q}_{0,\ep} 
	\to \vc{q}_0 \ \mbox{weakly in}\ L^1(R^3; R^3), \br  & 		
	\intRd{ \frac{1}{2} \frac{ |\vc{q}_{0, \ep}|^2 }{\vr_{0,\ep}}  } + 
	\int_{R^3 \setminus B_{\ep, 0}} P(\vr_{0, \ep})  \dx 
	\to  
	\intRd{ \frac{1}{2} \frac{ |\vc{q}_{0}|^2  }{\vr_0} } + 
	\intRd{ P(\vr_{0}) }  
	\label{P11}
	\end{align}	
as $\ep \to 0$.

		Then there is a subsequence (not relabelled) such that
		\begin{align}
		\vr_{\ep,f} &\to \vr \ \mbox{in}\ C_{\rm weak}([0,T]; L^\gamma (R^3)) \ \mbox{and in}\ L^1_{\rm loc}([0,T] \times R^3), \br
		\vue &\to \vu \ \mbox{weakly in}\ L^2(0,T; D^{1,2} (R^3; R^3)),
		\nonumber
		\end{align}
		where $(\vr, \vu)$ is a weak solution to the Navier--Stokes system \eqref{P1}--\eqref{P4} with the initial data $\vr_0$, $\vc{q}_0$.
		
		\end{Theorem}

	\end{mdframed}

\begin{Remark} \label{Rr1}
	
	Note that we may consider $\underline{\beta} = \Ov{\beta} = 0$ in hypothesis \eqref{P10} as soon as $\gamma > 3$.
	
	\end{Remark}

\begin{Remark} \label{Rr2}
	
	Here and hereafter, the symbol $a \aleq b$ means there is a positive constant $C$ such that $a \leq C b$.
	
\end{Remark}

The rest of the paper is devoted to the proof of Theorem \ref{PT1}. The leading idea is to use the test functions $\bfphi$ in the momentum equation \eqref{P7} that are constant (spatially homogeneous) 
on a neighborhood of the rigid body, in particular they satisfy \eqref{P8}. More specifically, the momentum balance yields that integral identity
\begin{align}
	\int_0^T &\intRd{ \Big[ \vre \vue \cdot \partial_t \bfphi + \vre \vue \otimes \vue : \Grad \bfphi + p(\vre) \Div \bfphi \Big]  } \dt
	\br &= \int_0^T \intRd{ \mathbb{S}(\Grad \vue) : \Grad \bfphi } \dt - \intRd{ \vc{q}_{0, \ep} \cdot \bfphi (0, \cdot) }
	\label{P13}
\end{align}
holds for any $\bfphi \in C^1_c([0,T) \times R^3; R^3)$ such that 
\begin{equation} \label{P14}
	\bfphi (t,x) = \avintBe{ \bfphi(t,\cdot) } \equiv \frac{1}{|B_{\ep,t}|} \int_{B_{\ep,t}} \bfphi (t,\cdot) \ \dx 
	\ \mbox{for any} \ x \ \mbox{in an open neighborhood of}\ B_{\ep, t}.
		\end{equation}
Using a simple density argument, it is easy to check that validity of \eqref{P13} can be extended to a larger class of test functions, namely 
$\bfphi \in W^{1,\infty}_c ([0,T) \times R^3; R^3)$, 
\begin{equation} \label{P14bis}
	\bfphi (t,x) = \avintBe{ \bfphi(t,\cdot) } \equiv \frac{1}{|B_{\ep,t}|} \int_{B_{\ep,t}} \bfphi (t,\cdot) \ \dx 
	\ \mbox{for any} \ x \in B_{\ep, t}.
\end{equation}
\section{Uniform bounds, weak convergence}
\label{U}

\subsection{Uniform bounds}

We start with uniform bounds that follow immediately from hypothesis \eqref{P11} and the energy inequality \eqref{P9}, namely
\begin{align}
{\rm ess} \sup_{t \in (0,T)} \|  \vr_{\ep,f}   \|_{L^1 \cap L^\gamma
	(R^3)} &\aleq 1,\  \label{U1a} \\ 
{\rm ess} \sup_{t \in (0,T)} \| \vre |\vue|^2 \|_{L^1(R^3)} &\aleq 1, \label{U1c} \\
{\rm ess} \sup_{t \in (0,T)} \| \vr_{\ep, f}\vue \|_{L^1 \cap L^{\frac{2 \gamma}{\gamma + 1}}(R^3; R^3)} &\aleq 1 \label{U1b}\\
\| \Ds \vue \|_{L^2(0,T; L^2 (R^3; R^{3 \times 3})) } &\aleq 1.
		\label{U1e}
		\end{align}
In particular, boundedness of the kinetic energy together with hypothesis 
\eqref{P10} yield the following estimate on the velocity of the rigid body 
\begin{equation} \label{U2}
\vr_{\ep, B} \ep^3 	|\vc{Y}_\ep (t) |^2 \aleq 1 \ \Rightarrow\ 
|\vc{Y}_\ep (t) | \aleq \ep^{\frac{1}{2}(\underline{\beta} - 3)}
 ,\ \vc{Y}_\ep = \frac{\D }{\dt} \vc{h}_\ep (t),\ t \in (0,T).
\end{equation}	

Finally, we deduce \eqref{U1e}
\begin{equation} \label{U1f}
\| \vue \|_{L^2(0,T; D^{1,2}(R^3; R^3))} \aleq 1\ \Rightarrow \ \| \vue \|_{L^2(0,T; L^6(R^3; R^3)) } \aleq 1.	
	\end{equation}
	
\subsection{ Convergence in continuity equation}
	
In view of the uniform bounds obtained in the preceding section, we deduce the existence of suitable subsequences satisfying
\begin{align} 
	\vr_{\ep, f} &\to \vr \ \mbox{in}\ C_{\rm weak}(0,T; L^\gamma (R^3)),  \br
	\vue &\to \vu \ \mbox{weakly in}\ L^2(0,T; D^{1,2} (R^d; R^d)), \br 
	\vr_{\ep,f} \vue &\to \vr \vu \ \mbox{weakly-(*) in}\ L^\infty(0,T; L^{\frac{2 \gamma}{\gamma + 1}} (R^3; R^3)),  
	\label{U3}
	\end{align}	
where we have used the fact that $(\vr_{\ep,f}, \vue)$ satisfy the equation of continuity \eqref{P5}, cf. \cite[Lemma 3.2]{EF64}. Now, it is easy to perform the limit 
in the equation of continuity \eqref{P5} to conclude
\begin{equation} \label{U4}
	\int_0^T \intRd{ \Big[ \vr \partial_t \varphi + \vr \vu \cdot \Grad \varphi \Big] } \dt = - \intRd{ \vr_{0} \varphi (0, \cdot) }
\end{equation}
for any $\varphi \in C^1_c([0,T) \times R^3)$.

Next, by virtue of hypothesis \eqref{P10}, 
\begin{equation} \label{U5}
\vre = \vr_{\ep,f} + \vr_{\ep,B} \mathds{1}_{B_{\ep}},\ \mbox{where}\ \vr_{\ep,B} \mathds{1}_{B_{\ep}} \to 0 \ \mbox{in}\ L^\infty(0,T; L^\Gamma(R^d)) 
\ \mbox{for some}\ \Gamma > \frac{3}{2}.
\end{equation}
In particular, the Young measure generated by $(\vr_{\ep,f})_{\ep > 0}$ coincides with that one generated by $(\vr_{\ep})_{\ep > 0}$. In particular, we may let $\ep \to 0$ in the renormalized equation of 
continuity \eqref{P6} obtaining
\begin{equation} \label{U6}
	\int_0^T \intRd{ \Big[ \Ov{b(\vr)} \partial_t \varphi + \Ov{b(\vr)} \vu \cdot \Grad \varphi + \Ov{( b(\vr) - b'(\vr) \vr )  \Div \vu } \varphi \Big] } \dt = - \intRd{ b(\vr_{0}) \varphi (0, \cdot) }
\end{equation}
for any $\varphi \in C^1_c([0,T) \times R^3)$ and any $b \in C^1[0, \infty)$, $b' \in C_c[0, \infty)$. Here and hereafter, the symbol $\Ov{b(\vr)}$ denotes the weak limit of the compositions 
$(b(\vre))_{\ep > 0}$ or, equivalently, $( b(\vr_{\ep,f}) )_{\ep > 0}$.

Finally, by the same token,
\[
\| \vr_{\ep,B} \mathds{1}_{B_{\ep, \tau}} \vue (\tau, \cdot) \|_{L^{\frac{2 \Gamma}{\Gamma + 1}}(R^3; R^3)} \leq \| \sqrt{\vr_{\ep, B}} \|_{L^{2 \Gamma}(B_{\ep, \tau})} \| \sqrt{\vre} \vue (\tau, \cdot) \|_{L^2(R^3; R^3)} \to 0 
\ \mbox{uniformly for}\ \tau \in (0,T);
\]
whence 
\begin{equation} \label{U7}
	\vr_{\ep,B}\mathds{1}_{B_{\ep, \tau}} \vue (\tau, \cdot) \to 0 \ \mbox{in}\ L^{\frac{2 \Gamma}{\Gamma + 1}}(R^3; R^3) \ \mbox{uniformly in}\ \tau \in (0,T).
\end{equation}	
Combining \eqref{U3}, \eqref{U5}, \eqref{U7} we may infer that
\begin{align} 
	\vre &\to \vr \ \mbox{in}\ C_{\rm weak}([0,T]; L^\gamma(R^3)) + L^\infty(0,T; L^\Gamma(R^3)), \label{U8a} \\ 
	\vre \vue &\to \vr \vu \ \mbox{in}\ L^{\infty}(0,T; L^{\frac{2 \gamma}{\gamma + 1}}(R^3; R^3))- \mbox{weak-(*)} + 
	L^\infty(0, T; L^{\frac{2 \Gamma}{\Gamma + 1}}(R^3; R^3)). 
\label{U8b}
	\end{align}
Moreover, 
\begin{align} 
	\| \vre \|_{L^\infty(0,T; L^1(R^3))} &\aleq 1 , \br 
	\| \vre \vue \|_{L^\infty(0,T; L^1(R^3; R^3))} &\aleq 1  .
	\label{U9}
\end{align}

\section{Restriction operators} 
\label{T}

In order to complete the proof of Theorem \ref{PT1} we have to address the following issues:
\begin{itemize} 
	\item convergence of the convective term $\vr_{\ep, f} \vue \otimes \vue$ to its counterpart $\vr \vu \otimes \vu$;
    \item uniform estimates and convergence of the pressure $p(\vr_{\ep,f})$; 
    \item limit passage $\ep \to 0$ in the momentum balance \eqref{P7}.
\end{itemize}

\noindent
To this end, we need a suitable restriction operator to accommodate the test functions in the class \eqref{P14}. As the result is of independent interest, we consider 
a general $d-$dimensional space, $d \geq 2$.

\subsection{ Construction of Restriction operator I}

Consider a function 
\begin{align} 
H &\in C^\infty(R), \ 0 \leq H(Z) \leq 1,\ H'(Z) = H'(1- Z) \ \mbox{for all}\ Z \in R, \br 
H(Z) &= 0 \ \mbox{for} \ - \infty < Z \leq \frac{1}{4},\ 
H(Z) = 1 \ \mbox{for}\ \frac{3}{4} \leq Z < \infty
\label{T1}
\end{align}

For $\varphi \in L^1(R^d)$ and $\vc{h} \in R^d$, we consider $E_\ep (\vc{h} )$,
\begin{align}
E_\ep (\vc{h})[\varphi] (x) = 
\frac{1}{|B_\ep (\vc{h}) |} \int_{B_\ep(\vc{h})} \varphi \ \dz \ H \left( 2 - \frac{|x - \vc{h}|}{\ep} \right)
+ \varphi(x) H \left( \frac{ | x - \vc{h} |}{\ep} - 1 \right).
\label{T2} 
\end{align}
where $B_\ep (\vc{h})$ denotes the ball centred at $\vc{h}$ with the radius $\ep > 0$.

The following properties are easy to check: 
\begin{itemize}
	\item 
The function $E_\ep (\vc{h})$ is constant on a small neighbourhood of the ball $B_\ep (\vc{h})$. Specifically, 
\begin{align}
	\frac{|x - \vc{h} |}{\ep} \leq 1 + \frac{1}{4} \ &\Rightarrow \ 
\frac{|x - \vc{h} |}{\ep} - 1 \leq \frac{1}{4} ,\ 
2 - \frac{|x - \vc{h} |}{\ep} \geq \frac{3}{4} 	\br &\Rightarrow 
E_\ep (\vc{h})[\varphi] =  \fint_{B_\ep(\vc{Y})} \varphi \ \dx.
\label{T3}		
	\end{align}
\item 
Similarly, 
\begin{align}
	\frac{|x - \vc{h} |}{\ep} \geq 1 + \frac{3}{4} \ &\Rightarrow \ 
	\frac{|x - \vc{h} |}{\ep} - 1 \geq \frac{3}{4} ,\ 
	2 - \frac{|x - \vc{h} |}{\ep} \leq \frac{1}{4} 	\br &\Rightarrow 
	E_\ep (\vc{h} )[\varphi] =  \varphi , 
	\label{T4}		
\end{align}
meaning $E_\ep (\vc{h})[\varphi]$ coincides with $\varphi$ on an open neighbourhood of the set 
$R^d \setminus B_{2 \ep}(\vc{h})$. 

\item 
\begin{equation} \label{T5}
{\rm supp} [ E_\ep (\vc{h})[\varphi] ] \subset \mathcal{U}_{3\ep} ({\rm supp}[\varphi])
\end{equation}
Indeed, if 
\[
{\rm dist}[ \vc{h}; {\rm supp}[\varphi]] \geq  \ep, 
\]
then 
\[
E_\ep (\vc{h})[\varphi] (x) = \varphi(x) H \left( \frac{ | x - \vc{h} |}{\ep} - 1 \right).  
\]
If
\[ 
{\rm dist}[ \vc{h}; {\rm supp}[\varphi]] < \ep,
\]
then, in accordance with \eqref{T4}, 
\[
E_\ep (\vc{h})[\varphi] (x) = \varphi (x) = 0 
\ \mbox{for a.a.}\ x \in R^d, {\rm dist}[x, {\rm supp}[\varphi]] > 3 \ep.
\]
\item 

In particular, it follows from \eqref{T5} that if $\varphi$ is compactly supported in an open set 
$\Omega \subset R^d$, then so is $E_\ep (\vc{h})[\varphi]$ provided $\ep > 0$ is small enough.

	\end{itemize}

Finally, by virtue of Jensen's inequality, 
\[
\left| \frac{1}{|B_\ep (\vc{h}) |} \int_{B_\ep(\vc{h})} \varphi \ \dz \right|^p \leq 
\frac{1}{|B_\ep (\vc{h}) |} \int_{B_\ep(\vc{h})} |\varphi|^p \ \dz,\ 1 \leq p < \infty.
\]
Consequently, we deduce 
\begin{equation} \label{lpe}
E_\ep (\vc{h}) [\varphi] = \varphi + \mathds{1}_{B_{\frac{7}{4} \ep} (\vc{h} ) } e^{0}_{\ep, \vc{h}}, \ \| e^{0}_{\ep, \vc{h}} \|_{L^p(R^d)} \aleq \| \varphi \|_{L^p (B_{\frac{7}{4} \ep} (\vc{h} ) ) },\ 
1 \leq p \leq \infty. 	
	\end{equation}

Summarizing we conclude that for any $\bfphi \in C^1_c([0,T) \times \Omega; R^d)$, the function 
$E_\ep (\vc{h}_\ep (\tau) )[\bfphi (\tau, \cdot)]$, $\tau \in [0,T]$ is an admissible test function in the momentum equation 
\eqref{P13}. Below, we derive the necessary error estimates on the spatial and time derivatives in Sobolev norms.

\subsubsection{Spatial derivatives}

Given $\vc{h} \in R^d$, the spatial derivatives of $E_{\ep}(\vc{h})$ can be computed directly using 
formula \eqref{T2}:
\begin{align} 
\Grad E_\ep (\vc{h})[\varphi] (x) &= \Grad \varphi (x) H \left( \frac{ | x - \vc{h} |}{\ep} - 1 \right) \br 
&+ \left( \varphi (x) - \frac{1}{|B_\ep (\vc{h}) |} \int_{B_\ep(\vc{h})} \varphi \ \dz \right) 
H' \left( \frac{ | x - \vc{h} |}{\ep} - 1 \right) \frac{1}{\ep} \frac{ x - \vc{h} }{|x - \vc{h}|},
\label{T6}
\end{align} 
where we have used that 
\[
H' \left( \frac{ | x - \vc{h} |}{\ep} - 1 \right) = 
H' \left( 2 - \frac{| x - \vc{h} |}{\ep}  \right).
\]

\subsubsection{Uniform bounds}

Seeing that 
\[
H' \left( \frac{ | x - \vc{h} |}{\ep} - 1 \right) \ne 0 
\ \Leftrightarrow \ \frac{5}{4} \ep \leq  | x - \vc{h} |  \leq \frac{7}{4} \ep
\]
we deduce 
\[
\left| \left( \varphi (x) - \frac{1}{|B_\ep (\vc{h}) |} \int_{B_\ep(\vc{h})} \varphi \ \dx \right) 
H' \left( \frac{ | x - \vc{h} |}{\ep} - 1 \right) \frac{1}{\ep} \frac{ x - \vc{h} }{|x - \vc{h}|} \right|
\aleq \| \Grad \varphi \|_{L^\infty(B_{\frac{7}{4} \ep}(\vc{h}) )} \mathds{1}_{B_{\frac{7}{4} \ep}(\vc{h})}
\]
Consequently, we deduce from \eqref{T6} the error estimates
\begin{align}
\Grad E_\ep (\vc{h})[\varphi]  = \Grad \varphi	+ e^{1}_{\ep, \vc{h}},\ 
|e^{1}_{\ep, \vc{h}}| \aleq \| \Grad \varphi \|_{L^\infty(B_{\frac{7}{4} \ep}(\vc{h}); R^d )} \mathds{1}_{B_{\frac{7}{4} \ep}(\vc{h})}
\label{T7}
\end{align}

\subsubsection{$L^p-$estimates on spatial derivatives}

Our goal is to show boundedness of the operator $E_\ep (\vc{h})$ in the Sobolev norms $W^{1,p}$. In view of 
formula \eqref{T6}, it is enough to control 
\[
 \left( \varphi (x) - \frac{1}{|B_\ep (\vc{h}) |} \int_{B_\ep(\vc{h})} \varphi \ \dz \right) 
H' \left( \frac{ | x - \vc{h} |}{\ep} - 1 \right) \frac{1}{\ep} \frac{ x - \vc{h} }{|x - \vc{h}|}
\]
on the annulus 
\[
\frac{5}{4} \ep \leq |x - \vc{h}| \leq \frac{7}{4} \ep
\]
in terms of the $L^p-$norm of $\Grad \bfphi$ on the same set. Without loss of generality, we may assume $\vc{h} = 0$. Thus our goal is to show the bound 
\begin{equation} \label{TT1}
\frac{1}{\ep} \left\|  \left( \varphi  - \frac{1}{|B_\ep  |} \int_{B_\ep} \varphi \ \dz \right) \right\|
_{L^p (\frac{5}{4} \ep \leq |x| \leq \frac{7}{4}\ep)} \aleq  
\| \Grad \varphi \|_{{L^p (B_{\frac{7}{4} \ep} }; R^{d \times d}) }	.
	\end{equation}
This is equivalent, after rescaling to the estimate 
\[
	 \left\|  \left( \varphi  - \frac{1}{|B_1  |} \int_{B_1} \varphi \ \dz \right) \right\|
	_{L^p (\frac{5}{4}  \leq |x| \leq \frac{7}{4})} \aleq  
	\| \Grad \varphi \|_{{L^p (B_{\frac{7}{4} } }; R^{d \times d}) }, 
\]
which, in turn, follows from Poincar\' e inequality 
\[
\| \varphi \|_{L^p(B_{\frac{7}{4}})} \aleq \| \Grad \varphi \|_{L^p(B_{\frac{7}{4}})} + 
\left| \int_{B_1} \varphi \ \D z \right|.
\]
Thus, together with \eqref{T7}, we conclude
\begin{equation} \label{TT7}
\Grad E_\ep (\vc{h})[\varphi]  = \Grad \varphi	+ \mathds{1}_{B_{\frac{7}{4} \ep}(\vc{h})} e_{\ep, \vc{h}}, \ \| e_{\ep, \vc{h}} \|_{L^p(R^{d})} \aleq \| \Grad \varphi \|_{{L^p (B_{\frac{7}{4} \ep} }; R^{d \times d}) } \ 1 \leq p \leq \infty.	
	\end{equation}

\subsubsection{Derivative with respect to the parameter $\vc{h}$}

Similarly to the preceding part, we compute
\begin{align} 
	\nabla_{\vc{h}} E_\ep (\vc{h})[\varphi] (x) &= \nabla_{\vc{h}} \left( \frac{1}{|B_\ep (\vc{h}) |} \int_{B_\ep(\vc{h})} \varphi \ \dz \right) H \left( 2 - \frac{|x - \vc{h}|}{\ep} \right) \br 
	&- \left( \varphi (x) - \frac{1}{|B_\ep (\vc{h}) |} \int_{B_\ep(\vc{h})} \varphi \ \dz \right) 
	H' \left( \frac{ | x - \vc{h} |}{\ep} - 1 \right) \frac{1}{\ep} \frac{ x - \vc{h} }{|x - \vc{h}|},
	\label{T8}
\end{align} 
where, furthermore, 
\[
\nabla_{\vc{h}} \left( \frac{1}{|B_\ep (\vc{h}) |} \int_{B_\ep(\vc{h})} \varphi \ \dz \right) = 
\frac{1}{|B_\ep (\vc{h}) |}\int_{B_\ep(\vc{h})} \Grad \varphi \ \dz.
\]
We therefore obtain 
\begin{align} 
	\nabla_{\vc{h}} E_\ep (\vc{h})[\varphi] (x) &=  \frac{1}{|B_\ep (\vc{h}) |} \int_{B_\ep(\vc{h})} \Grad \varphi \ \dz  H \left( 2 - \frac{|x - \vc{h}|}{\ep} \right) \br 
	&- \left( \varphi (x) - \frac{1}{|B_\ep (\vc{h}) |} \int_{B_\ep(\vc{h})} \varphi \ \dz \right) 
	H' \left( \frac{ | x - \vc{h} |}{\ep} - 1 \right) \frac{1}{\ep} \frac{ x - \vc{h} }{|x - \vc{h}|},
	\label{T9}
\end{align} 
which can be also written as a commutator
\begin{align} 
	\nabla_{\vc{h}} E_\ep (\vc{h})[\varphi] (x) &=  \frac{1}{|B_\ep (\vc{h}) |} \int_{B_\ep(\vc{h})} \Grad \varphi \ \dz  H \left( 2 - \frac{|x - \vc{h}|}{\ep} \right) + H \left( \frac{ | x - \vc{h} |}{\ep} - 1 \right) \Grad \varphi
	\br 
	&- \Grad E_\ep (\vc{h}) [\varphi] = E_\ep (\vc{h})[\Grad \varphi ]  - \Grad E_\ep (\vc{h}) [\varphi].
	\label{T10}
\end{align}

\subsubsection{Estimates on the time derivative}

The time derivative of the restriction operator $E_{\ep}(\vc{h}(\tau))[\bfphi (\tau, \cdot)]$ can be computed by using formula \eqref{T10}: 
\begin{align}
\partial_t &\left( E_\ep (\vc{h} (t))[ \varphi (t, \cdot) ] \right) = E_\ep (\vc{h} (t))[ \partial_t \varphi (t, \cdot) ]
+ \nabla_{\vc{h}} E_\ep (\vc{h} (t))[ \varphi (t, \cdot) ] \cdot \frac{\D }{\dt} \vc{h} (t) \br
&= E_\ep (\vc{h} (t))[ \partial_t \varphi (t, \cdot) ] + E_\ep (\vc{h}(t))[\Grad \varphi(t, \cdot)  ]\cdot \vc{Y} (t) 
 - \Grad E_\ep (\vc{h}(t)) [\varphi (t, \cdot)] \cdot \vc{Y} (t),
\label{T11}
\end{align}
where
	\[
\vc{Y} = \frac{\D }{\dt} \vc{h}.
\]

We conclude this section by summarizing the basic properties of the restriction operator $E_\ep (\vc{h})$. 

\begin{Proposition} \label{EP1a}
	The operator $E_\ep (\vc{h})[\varphi]$ is well defined for $\varphi \in L^1_{\rm loc}(R^d)$. The following holds true:
	\begin{itemize}
		\item
		\begin{equation} \label{Ta1}
		 E_\ep (\vc{h})[\varphi] = \left\{ \begin{array}{l} \frac{1}{|B_{\ep}(\vc{h})|} 
		 	\int_{B_{\ep}(\vc{h})} \varphi \ \dx \ \mbox{if}\ |x - \vc{h}| < \ep  \\ 
		 	\\ 
		 	\varphi \ \mbox{if}\ |x - \vc{h}| > 2 \ep  \end{array}	\right. ;
			\end{equation}
		\item
		\begin{equation} \label{Ta2} 
	E_\ep (\vc{h}) [\varphi] = \varphi + \mathds{1}_{B_{\frac{7}{4} \ep} (\vc{h}) } e^0_{\ep, \vc{h}}, \ \| e^0_{\ep, \vc{h}} \|_{L^p(R^d)} \aleq \| \varphi \|_{L^p (B_{\frac{7}{4} \ep} (\vc{Y} ) ) },\ 
	1 \leq p \leq \infty; 	
	\end{equation}
\item
\begin{equation} \label{Ta3}
\Grad E_\ep (\vc{h})[\varphi]  = \Grad \varphi	+ \mathds{1}_{B_{\frac{7}{4} \ep}(\vc{h})} e^1_{\ep, \vc{h}}, \ \| e^1_{\ep, \vc{h}} \|_{L^p(R^{d})} \aleq \| \Grad \varphi \|_{{L^p (B_{\frac{7}{4} \ep}(\vc{h}) }; R^{d \times d}) } \ 1 \leq p \leq \infty;
\end{equation}
		\item if $\vc{h} = \vc{h}(t)$ is Lipschitz and $\varphi \in W^{1,1}_{\rm loc}( (0,T) \times R^d)$, then 
	\begin{align}
		\partial_t &\left( E_\ep (\vc{h} (t))[ \varphi (t, \cdot) ] \right) \br
		&= E_\ep (\vc{h} (t))[ \partial_t \varphi (t, \cdot) ] + E_\ep (\vc{h}(t))[\Grad \varphi(t, \cdot) ] \cdot \vc{Y} (t) 
		- \Grad E_\ep (\vc{Y}(t)) [\varphi (t, \cdot)] \cdot \vc{Y} (t)
		\label{Ta4}
	\end{align}	
for a.a. $t \in (0,T)$,
		where 
		\[
		\vc{Y} = \frac{\D }{\dt} \vc{h}.
		\]
		
		\end{itemize}

	\end{Proposition}

\section{Restriction operator revisited}
\label{E}

The drawback of the restriction operator $E_\ep (\vc{h}(t))$ is that it does not preserve the divergence of a vector valued function. To remedy this, we introduce a new vector--valued restriction operator $\vc{R}_\ep (\vc{h})$ acting on vector valued functions.

\subsection{Basic structure}

We start by introducing {\em shift} operator 
\[
S_{\vc{h}} [f] (x) = f ( \vc{h} + x ).
\]
Setting 
\[
E_{\ep} = E_{\ep} (\vc{0}) 
\]
for the restriction operator introduced in the previous section, we check easily the relation 
\[
E_{\ep} (\vc{h}) [ \varphi ] = S_{- \vc{h}} E_{\ep}\Big[ S_\vc{h} [\varphi] \Big].
\]
We compute 
\begin{equation} \label{E1}
\nabla_{\vc{h}} S_{\vc{h}} [f] (x) = \Grad f (\vc{h} + x ) = S_{\vc{h}} [ \Grad f ],
\end{equation}
which, in particular, yields the commutator formula \eqref{T10}. 

\subsection{Bogovskii operator} 

We use a particular version of Bogovskii operator constructed by Diening, R\accent23u\v zi\v cka, Schumacher \cite{DieRuzSch}. 
The operator $\mathcal{B}_{2\ep,  \ep}$ is a branch of the inverse of the divergence operator defined on the 
annulus, 
\[
B_{2\ep} \setminus B_{\ep}. 
\]
The operator enjoys the following properties: 
\begin{itemize}
	\item 
	\begin{equation} \label{E2}
	\mathcal{B}_{2 \ep, \ep} : L^p_0 (B_{2\ep} \setminus B_{\ep}) \to W^{1,p}_0 
(B_{2\ep} \setminus B_{\ep}),
\end{equation}
$L^p_0$ denoting the space of $L^p$ functions with zero mean,\ 
\begin{equation} \label{E3}
	\| \Grad \mathcal{B}_{2 \ep, \ep} [f] \|_{L^p (B_{2\ep} \setminus B_{\ep}; R^{d \times d} )} \aleq 
	\| f \|_{L^p (B_{2\ep} \setminus B_{\ep} ) },		
		\end{equation}
for any $1 < p < \infty$, where the embedding constant is independent of $\ep$, 
\begin{equation} \label{E4}
	\Div \mathcal{B}_{2 \ep, \ep}[f] = f \ \mbox{in}\ B_{2\ep} \setminus B_{\ep}. 	
	\end{equation}
\item If, in addition, $f = \Div \vc{g}\in L^q_0 (B_{2\ep} \setminus B_{\ep})$, where $\vc{g} \in W^{1,q} (B_{2\ep} \setminus B_{\ep})$,\
$\vc{g} \cdot \vc{n}|_{\partial (B_{2\ep} \setminus B_{\ep})} = 0$,
then 
\begin{equation} \label{E5}
\| \mathcal{B}_{2 \ep, \ep}[\Div \vc{g} ] \|_{L^q (	B_{2\ep} \setminus B_{\ep} ;R^d) } \aleq 
\|  \vc{g}  \|_{L^q (	B_{2\ep} \setminus B_{\ep} ;R^d) } 
	\end{equation}
for $1 < q < \infty$, where the embedding constant is independent of $\ep$.
	\end{itemize}

The operator $\mathcal{B}_{2 \ep, \ep}$ was constructed by Diening et al \cite{DieRuzSch}. The remarkable property that its norms are independent of $\ep$ follow from the fact that $B_{2\ep} \setminus B_{\ep}$ are John domains 
uniformly in $\ep$, see Diening et al \cite[Theorem 5.2]{DieRuzSch}, Lu and Schwarzacher \cite[Theorem 1.1]{LuSchw}.

\subsection{ Construction of Restriction operator II}
We define the operator 
\begin{equation} \label{E6}
	\vc{R}_\ep [\bfphi ] = E_\ep [\bfphi] + \mathcal{B}_{2 \ep, \ep} \Big[ \left(\Div\bfphi - 
	\Div E_\ep [\bfphi]\right)|_{B_{2\ep} \setminus B_{\ep}} \Big]
	\end{equation}
{\it A priori} the operator is defined for $\bfphi \in W^{1,p}(R^d; R^d)$. For the definition to be correct, 
we have to verify that $\left(\Div\bfphi - 
	\Div E_\ep [\bfphi]\right)$ has zero mean over the annulus $B_{2\ep} \setminus B_{\ep}$. It is enough if we
consider the functions $\bfphi$ with the following properties: 
\begin{equation} \label{E7}
\bfphi \cdot \vc{n}|_{|x| = 2 \ep} = E_{\ep} [\bfphi ] \cdot \vc{n}|_{|x| = 2 \ep }, 
\end{equation}
and 
\begin{equation} \label{E8}
	\int_{|x| = \ep } \bfphi \cdot \vc{n} \ \D \sigma = 	\int_{|x| = \ep } E_{\ep} [\bfphi ] \cdot \vc{n} \ \D \sigma.
\end{equation}
On the one hand, 
equality \eqref{E7} obviously holds as $\bfphi = E_{\ep}[ \bfphi]$ if $|x| = 2 \ep$. On the other hand, 
we have
\[ 
\int_{|x| = \ep }  E_{\ep} [\bfphi ] \cdot \vc{n} \ \D \sigma = 0
\]
as $E_{\ep} [\bfphi ]|_{B_\ep}$ is constant. Thus for \eqref{E8} to hold, it is enough to assume 
\begin{equation} \label{E9}
	\Div \bfphi|_{B_\ep} = 0.
	\end{equation}

Finally, we set 
\begin{equation} \label{E10}
	\vc{R}_\ep (\vc{h} )[\bfphi] = S_{- \vc{h}} \vc{R}_\ep \Big[ S_{\vc{h}} [\bfphi] \Big].
	\end{equation}
Summarizing the previous discussion, we get.

\begin{Proposition}[{\bf Continuity in $L^p$ spaces}] \label{EP1}
	
	The operator $\vc{R}_\ep (\vc{h})$ is well defined for any function $\bfphi \in W^{1,p}(R^d; R^d)$ 
satisfying 
\begin{equation} \label{E11}
	\Div\bfphi = 0 \ \mbox{for all}\ x, \ |x - \vc{h}| < \ep .
	\end{equation}
Moreover, 
\begin{itemize} 
	\item 
\begin{equation} \label{E12}
\vc{R}_\ep (\vc{h})[\bfphi] = \left\{ \begin{array}{l} \frac{1}{|B_{\ep}(\vc{h})|} 
\int_{B_{\ep}(\vc{h})} \bfphi \ \dx \ \mbox{if}\ |x - \vc{h}| < \ep,  \\ 
\\ 
\bfphi \ \mbox{if}\ |x - \vc{h}| > 2 \ep;  \end{array}	\right. 
	\end{equation} 
\item 
\begin{equation} \label{E13}
	\Div \vc{R}_\ep (\vc{h})[\bfphi] = \Div \bfphi ;
	\end{equation}
\item 
\begin{equation} \label{E14}
\| \vc{R}_\ep (\vc{h})[\bfphi] \|_{W^{1,p}(R^d; R^d)} \aleq 
\| \bfphi\|_{W^{1,p}(R^d; R^d)}
\end{equation} 
for any $1 < p < \infty$ independently of $\ep > 0$.
	
	\end{itemize}

\end{Proposition}

\subsection{Estimates in the negative norm}

In order to estimate time derivatives, we need to find bounds on $\vc{R}_\ep$ provided the argument is 
in the form 
\begin{equation} \label{E15}
\bfphi = \mathcal{B}[ \Div \vc{g} ] , 
\end{equation}
where $\mathcal{B}$ is some right inverse of the divergence operator, $\Div \circ \mathcal{B} = {\rm Id}$, 
such that 
\[
\mathcal{B} \circ \Div \ \mbox{is bounded in}\ L^q(R^d; R^d) \ \mbox{for any}\ 1 < q < \infty,
\]
and
\[
\vc{g}|_{B_\ep (\vc{h}) } = 0.
\]
If this is the case, then 
\[
\Div \bfphi = \Div \mathcal{B}[ \Div \vc{g} ] = \Div \vc{g} = 0 \ \mbox{on}\ B_\ep (\vc{h}) 
\]
so the operator $\vc{R}_\ep (\vc{h})[ \bfphi]$ is well defined. Below, we consider $\mathcal{B} = \Grad \Del^{-1}$, however, $\mathcal{B}$ can be also the standard Bogovskii operator on some domain $\Omega \subset R^d$.

Without loss of generality, we may assume 
\begin{equation} \label{E16}
\vc{h} = 0,\ \vc{R}_\ep (\vc{h}) = \vc{R}_\ep ,\ 
\vc{g}|_{B_\ep} = \Div \bfphi|_{B_\ep} = 0.
\end{equation}

Our goal is to obtain $L^q$ estimates on $\vc{R}_\ep [\bfphi]$ in terms of the $L^q-$norm of $\vc{g}$. 
As $E_\ep$ is bounded as an operator on $L^q$, it is enough to check boundedness of the term 
\begin{align}
&\mathcal{B}_{2 \ep, \ep} \left[ \Div \mathcal{B}[ \Div \vc{g} ] - 
\Div E_\ep \Big[ \mathcal{B}\left[ \Div \vc{g} \right] \Big] \right] \br 
&=  \mathcal{B}_{2 \ep, \ep} \Big[ \Big( E_{\ep} \left[ \Div \mathcal{B}[ \Div \vc{g} ] \right]
- 
\Div E_\ep \Big[ \mathcal{B}\left[ \Div \vc{g} \right] \Big] \Big) + 
\Big( \Div \mathcal{B}[ \Div \vc{g} ] - \Big( E_{\ep} \left[ \Div \mathcal{B}[ \Div \vc{g} ] \right]\Big) 
\Big].
\label{E17}
\end{align}

We get 
\begin{align}
E_{\ep} &\left[ \Div \mathcal{B}[ \Div \vc{g} ] \right]
- 
\Div E_\ep \Big[ \mathcal{B}\left[ \Div \vc{g} \right] \Big]	\br
&=  E_{\ep} \left[  \Div \vc{g}  \right] - 
\Div \left( \mathcal{B} [\Div \vc{g} ] H \left( \frac{|x|}{\ep} - 1 \right) 
+ H \left( 2 - \frac{|x|}{\ep} \right) \frac{1}{|B_\ep |} \int_{B_\ep} \mathcal{B}[\Div \vc{g}] \dx   \right)
\br 
&= \Div \vc{g} H \left( \frac{|x|}{\ep} - 1 \right) 
- \Div \vc{g} H \left( \frac{|x|}{\ep} - 1 \right) - \frac{1}{\ep} \mathcal{B}[\Div \vc{g}] H'\left( \frac{|x|}{\ep} - 1 \right) \frac{x}{|x|} \br 
&+ \frac{1}{\ep} \frac{x}{|x|} H' \left( 2 - \frac{|x|}{\ep} \right) \frac{1}{|B_\ep|}\int_{B_\ep} \mathcal{B}[\Div \vc{g}] \dx \br 
&= \frac{1}{\ep} \frac{x}{|x|}  H'\left( \frac{|x|}{\ep} - 1 \right) 
\left( \frac{1}{|B_\ep|}\int_{B_\ep} \mathcal{B}[\Div \vc{g}] \dx  - \mathcal{B}[\Div \vc{g} ] \right). 
\label{E18}
\end{align}
Furthermore,  
\begin{align}
 \Div \mathcal{B}[ \Div \vc{g} ] &-  E_{\ep} \left[ \Div \mathcal{B}[ \Div \vc{g} ] \right]	 = 
 \Div \vc{g} -  E_{\ep} \left[  \Div \vc{g}  \right] \br 
& = \Div \vc{g} - \Div \vc{g} H \left( \frac{|x|}{\ep} - 1 \right) - H \left( 2 - \frac{|x|}{\ep} \right) 
 \frac{1}{|B_\ep |} \int_{B_\ep} \Div \vc{g} \dx \br
&=  \Div \vc{g} - \Div \vc{g} H \left( \frac{|x|}{\ep} - 1 \right) \br
&= \Div \left( \vc{g} - \vc{g} H \left( \frac{|x|}{\ep} - 1 \right) \right) + 
\frac{1}{\ep} \vc{g} H' \left (\frac{|x|}{\ep} - 1 \right) \frac{x}{|x|} 
	\label{E19}
	\end{align}
Summing up the previous relations, we conclude 
\begin{align}
\mathcal{B}_{2 \ep, \ep} &\left[ \Div \mathcal{B}[ \Div \vc{g} ] - 
\Div E_\ep \Big[ \mathcal{B}\left[ \Div \vc{g} \right] \Big] \right] \br 
&= \mathcal{B}_{2 \ep, \ep} \left[ \frac{1}{\ep} \frac{x}{|x|}  H'\left( \frac{|x|}{\ep} - 1 \right) 
\left( \frac{1}{|B_\ep|}\int_{B_\ep} \mathcal{B}[\Div \vc{g}] \dx  - \mathcal{B}[\Div \vc{g} ] \right) \right. \br 
&+ \left. 	\frac{1}{\ep} \vc{g} H' \left (\frac{|x|}{\ep} - 1 \right) \frac{x}{|x|} \right] \br 
&+ \mathcal{B}_{2 \ep, \ep} \left[ \Div \left( \vc{g} - \vc{g} H \left( \frac{|x|}{\ep} - 1 \right) \right) \right].
	\label{E20}
	\end{align}

Seeing that $\vc{g}|_{B_\ep} = 0$ we may use the ``negative'' estimates \eqref{E5} to deduce 
\begin{equation} \label{E21} 
\left\|  \mathcal{B}_{2 \ep, \ep} \left[ \Div \left( \vc{g} - \vc{g} H \left( \frac{|x|}{\ep} - 1 \right) \right) \right]	\right\|_{L^p (B_{2\ep} \setminus B_\ep ; R^d) } \aleq 
\left\| \vc{g} \right\|_{L^p(B_{2\ep} \setminus B_\ep; R^d ) },\ 1 < p < \infty. 
	\end{equation}
Finally, by means of the $L^p-$bounds \eqref{E3}, 
\begin{align} 
&\left\| \Grad \mathcal{B}_{2 \ep, \ep} \left[ \frac{1}{\ep} \frac{x}{|x|}  H'\left( \frac{|x|}{\ep} - 1 \right) 
\left( \frac{1}{|B_\ep|}\int_{B_\ep} \mathcal{B}[\Div \vc{g}] \dx  - \mathcal{B}[\Div \vc{g} ] \right) \right. \right.	\br
&\quad +\left. \left.	\frac{1}{\ep} \vc{g} H' \left (\frac{|x|}{\ep} - 1 \right) \frac{x}{|x|} \right]
 \right\|_{L^p(B_{2\ep} \setminus B_\ep; R^{d \times d} ) } \br \quad& \aleq 
 \frac{1}{\ep} \left( \| \mathcal{B}[\Div \vc{g}] \|_{L^p(B_{2 \ep} \setminus B_\ep ; R^d)} + 
\| \vc{g}\|_{L^p(B_{2 \ep} \setminus B_\ep ; R^d)} \aleq \frac{1}{\ep} \| \vc{g}\|_{L^p(R^d ; R^d)}
 \right).
 \label{E22}
	\end{align}
Thus, by virtue of Poincar\` e inequality on $B_{2\ep} \setminus B_\ep$, 
\begin{align} 
	&\left\| \mathcal{B}_{2 \ep, \ep} \left[ \frac{1}{\ep} \frac{x}{|x|}  H'\left( \frac{|x|}{\ep} - 1 \right) 
	\left( \frac{1}{|B_\ep|}\int_{B_\ep} \mathcal{B}[\Div \vc{g}] \dx  - \mathcal{B}[\Div \vc{g} ] \right) \right. \right.	\br
	&\quad +\left. \left.	\frac{1}{\ep} \vc{g} H' \left (\frac{|x|}{\ep} - 1 \right) \frac{x}{|x|} \right]
	\right\|_{L^p(B_{2\ep} \setminus B_\ep; R^{d} ) }  \aleq 
	 \| \vc{g} \|_{L^p(R^d ; R^d)}.
	\label{E32}
\end{align}

We have obtained the following result.

\begin{Proposition} [{\bf Continuity in the negative space}] \label{EP2}
Let $\bfphi \in W^{1,p}_{\rm loc}(R^d; R^d)$ can be written in the form 
\[
\bfphi = \mathcal{B}[\Div \vc{g} ],\ \vc{g} \in L^q(R^d; R^d), \ \vc{g}|_{B_\ep (\vc{h})} = 
 0,
\]	
where $\Div \circ \mathcal{B} = {\rm Id}$ and $\mathcal{B} \circ \Div$ bounded on $L^q(R^d)$, $1 < q < \infty$.

Then 
\[ 
\left\| \vc{R}_\ep (\vc{h}) [ \mathcal{B}[ \Div \vc{g} ] ]  \right\|_{L^q(R^d; R^d)} \aleq 
\left\| \vc{g}  \right\|_{L^q(R^d; R^d)} ,\ 1 < q < \infty
\]
uniformly in $\ep$.
	\end{Proposition}

\begin{Remark} \label{ER1}
	The same result holds for general operators of the form 
	\[
	\bfphi = \mathcal{L} [\vc{g}] , 
	\]
	provided 
	\[
	\vc{g}|_{B_\ep (\vc{h})} = 0,\ 
	\Div  \mathcal{L} [\vc{g}] = \Div \vc{g}.
	\]
	In particular, we may consider 
	\[
	\bfphi = \Grad \mathcal{B}[r] \cdot \vc{V},\ \mbox{with} \ \vc{g} = r \vc{V}, \ 
	\varphi_i = \partial_{x_j} \mathcal{B}_i [r] \vc{V}_j,
	\]
	where $\vc{V}\in R^d$ is a constant vector. 
	\end{Remark}

\section{Pressure estimates}
\label{pe}

The well known problem connected with the compressible fluid flow is the lack of integrability of the pressure term $p(\vr)$ in the $x-$variable. 
If $\gamma > \frac{3}{2}$, the relevant estimates are obtained considering the quantity
\begin{equation} \label{pe1a}
\bfphi = \vc{R}_\ep (\vc{h}_\ep (t) ) [ \Grad \Del^{-1} b(\vr_{\ep,f}) ]
\end{equation}
as a test function in the momentum equation \eqref{P13}, where 
\[
b(r) \geq 0,\ b(r) = 0 \ \mbox{for all}\ 0 \leq r \leq 1,\ b(r) = r^\alpha \ \mbox{for}\ r \geq 2, \ \alpha\in (0,\gamma),
\]
and $\Del^{-1}$ denotes the inverse of the Laplace operator on $R^3$,  
\[
\Del^{-1}[v] = \mathcal{F}^{-1}_{\xi \to x} \left[ \frac{1}{|\xi|^2} \mathcal{F}_{x \to \xi} [v] \right],\ 
\mathcal{F} - \mbox{the Fourier transform.}
\]
Note carefully that 
\[
\Div [ \Grad \Del^{-1} b(\vr_{\ep,f})(t, \cdot) ] = b(\vr_{\ep,f})(t,\cdot) = 0 \ \mbox{on}\ B_{\ep,t}.
\]

In accordance with the uniform bounds \eqref{U1a},
\begin{equation} \label{pe1}
	{\rm ess} \sup_{t \in (0,T)}
\| b(\vr_{\ep, f})(t, \cdot) \|_{L^1 \cap L^{\frac{\gamma}{\alpha}} (R^3)}   \aleq 1.	
	\end{equation}
Moreover, evoking the standard elliptic estimates, we get
\begin{align}
\Grad [ \Grad \Del^{-1} b(\vr_{\ep,f}) ] \ &\mbox{bounded in}\ L^\infty(0,T; L^r(R^3; R^{3 \times 3})) \ \mbox{for any}\ 1 < r \leq \frac{\gamma}{\alpha}, \br 
[ \Grad \Del^{-1} b(\vr_{\ep,f}) ] \ &\mbox{bounded in}\  L^\infty(0,T; L^q(R^d; R^{d})) \ \mbox{for any}\ \frac{3}{2} < q \leq \infty	
	\label{pe2}
	\end{align}
provided $\frac{\gamma}{\alpha} > 3$.

\subsection{Equi--integrability of the pressure}

Using $\bfphi$ introduced in \eqref{pe1a} as a test function in the momentum balance \eqref{P13} we get 
\begin{align}
\int_0^T &\psi(t) \intRd{ p(\vr_{\ep,f}) b(\vr_{\ep,f}) } \dt \br &= 
\int_0^T \psi (t) \intRd{ p(\vr_{\ep,f}) \Div \vc{R}_\ep (\vc{h}_\ep (t) ) [ \Grad \Del^{-1} b(\vr_{\ep,f}) ] } \dt = \sum_{i=1}^5 I_{i,\ep},
 \label{pe4}
\end{align}		
for any $\psi \in C^1_c[0,T)$, $\psi(0) = 1$, where 
\begin{align}
I_{1,\ep} &= \int_0^T \psi \intRd{ \mathbb{S}(\Grad \vue) : \Grad \vc{R}_\ep (\vc{h}_\ep (t) ) [ \Grad \Del^{-1} b(\vr_{\ep,f}) ] } \dt, \br
I_{2,\ep} &= - \int_0^T \psi \intRd{ \vr_{\ep,f} \vue \otimes \vue : \Grad \vc{R}_\ep (\vc{h}_\ep (t) ) [ \Grad \Del^{-1} b(\vr_{\ep,f}) ] } \dt, \br
I_{3,\ep} &= - \intRd{ \vc{q}_{0, \ep} \cdot \vc{R}_\ep (\vc{Y}_\ep (0) ) [ \Grad \Del^{-1} b(\vr_{\ep,f}) (0, \cdot) ] }, \br
I_{4,\ep} &= \int_0^T \partial_t \psi \intRd{ \vr_{\ep} \vue \cdot  \vc{R}_\ep (\vc{h}_\ep (t) ) [ \Grad \Del^{-1} b(\vr_{\ep,f}) ] }, \br
I_{5,\ep} &= \int_0^T  \psi \intRd{ \vr_{\ep} \vue \cdot  \partial_t \vc{R}_\ep (\vc{h}_\ep (t) ) [ \Grad \Del^{-1} b(\vr_{\ep,f}) ] }. 
	\label{pe5}
 \end{align}
Our goal is to show that all integrals $I_{i,\ep}$, $i=1,\dots,5$ are bounded uniformly for $\ep \to 0$ as soon as $\alpha > 0$ is chosen small enough. Accordingly, relation 
\eqref{pe4} together with the bound \eqref{U1a}, yield equi--integrability of the pressure
\begin{equation} \label{equip}
	\int_0^T \intRd{ p(\vr_{\ep,f}) \vr_{\ep,f}^\alpha } \dt \aleq 1.
\end{equation}

\subsubsection{Viscosity and convective term}

It follows form \eqref{E14} and \eqref{pe2} that 
\[
\left( \Grad \vc{R}_\ep (\vc{h}_\ep (t) ) [ \Grad \Del^{-1} b(\vr_{\ep,f}) ] \right)_{\ep > 0} \ \mbox{is bounded in}\ L^\infty(0,T; L^r(R^{3 \times 3})),\ 
1 < r \leq \frac{\gamma}{\alpha}.
\]
In particular, the integral $I_{1,\ep}$ remains bounded uniformly for $\ep \to 0$. 

Similarly, in view of the energy estimates \eqref{U1c}, 
\[
\left( \vr_{\ep, f} \vue \otimes \vue \right)_{\ep > 0}\ \mbox{is bounded in} \ L^\infty(0,T; L^1(R^3, R^{3 \times 3})).
\]
Morever, as $\vue$ is bounded in $L^2(0,T; L^6(R^3; R^3))$ (see \eqref{U1f}), we get 
\[
\left( \vr_{\ep, f} \vue \otimes \vue \right)_{\ep > 0} 
\ \mbox{bounded in}\ L^1(0,T; L^s (R^3, R^{3 \times 3})),\ s > 1,\ \frac{1}{s} = \frac{1}{3} + \frac{1}{\gamma}.
\] 
Consequently, by interpolation, 
\begin{equation} \label{inter}
\left( \vr_{\ep, f} \vue \otimes \vue \right)_{\ep > 0}\ \mbox{is bounded in} 
\\ L^q(0,T; L^m (R^3, R^{3 \times 3}))\ \mbox{for some}\ q > 1, m > 1.
\end{equation}
In particular, $I_{2,\ep}$ remains bounded uniformly for $\ep \to 0$.

\subsubsection{Momentum}

In accordance with the bounds \eqref{U1b}, \eqref{U7}, we obtain
\[
(\vre \vue)_{\ep > 0} \ \mbox{is bounded in}\ L^\infty \left(0,T; \left( L^{\frac{2 \gamma}{\gamma + 1}} + L^{\frac{2 \Gamma}{\Gamma + 1}} \right) (R^3; R^3) \right).
\]
Moreover, the relation \eqref{pe2}$_2$ gives
\[
\left( \vc{R}_\ep (\vc{h}_\ep (t) ) [ \Grad \Del^{-1} b(\vr_{\ep,f}) ] \right)_{\ep > 0} \ \mbox{is bounded in}\ 
L^\infty(0,T; L^q(R^3;R^3)) \ \mbox{for any}\ \frac{3}{2} < q \leq \infty.
\]
Thus we conclude that $I_{\ep,4}$, and, similarly, $I_{3,\ep}$ remain bounded for $\ep \to 0$.

\subsubsection{Time derivative}

In order to evaluate the time derivative in $I_5$, we have an analogue of formula \eqref{T11} by using the relations \eqref{E1} and \eqref{E10}: 
\begin{align}
\partial_t \vc{R}_\ep (\vc{h}_\ep (t) ) &[ \Grad \Del^{-1} b(\vr_{\ep,f}) ] = 
\vc{R}_\ep (\vc{h}_\ep (t) ) [ \Grad \Del^{-1} \partial_t b(\vr_{\ep,f}) ]	\br &+ 
\vc{R}_\ep (\vc{h}_\ep (t) ) [ \Grad \Grad \Del^{-1} b(\vr_{\ep,f})] \cdot \vc{Y}_\ep   - 
\Grad \vc{R}_\ep (\vc{h}_\ep (t) ) [ \Grad \Del^{-1} b(\vr_{\ep,f}) ] \cdot \vc{Y}_\ep, 
	\label{pe6}
	\end{align}
where $\vr_{\ep,f}$ satisfies the renormalized equation of continuuity, 
\[
\Grad \Del^{-1} \partial_t b(\vr_{\ep,f}) = - \Grad \Del^{-1} \Div (b(\vr_{\ep,f}) \vue) + \Grad \Del^{-1} \left[ \Big( b(\vr_{\ep, f}) - b'(\vr_{\ep, f}) \vr_{\ep, f} \Big) \Div \vue \right].
\]

Now, in accordance with Proposition \ref{EP2} and Remark \ref{ER1}, we get 
\begin{align} 
\| \vc{R}_\ep (\vc{h}_\ep (t) ) [ \Grad \Grad \Del^{-1} b(\vr_{\ep,f})] \cdot \vc{Y}_\ep  	\|_{L^r(R^3)} \aleq |\vc{Y}_\ep|\ \| b(\vr_{\ep,f}) \|_{L^r(R^3)} 
\aleq |\vc{Y}_\ep|,\ 1 < r \leq \frac{\gamma}{\alpha}, \br 
\| \vc{R}_\ep (\vc{h}_\ep (t) ) [\Grad \Del^{-1} \Div (b(\vr_{\ep,f}) \vue)   ]	\|_{L^q(R^3)} \aleq \| b(\vr_{\ep,f}) \vue \|_{L^q(R^3)},\ 1 < q < \infty.
	\label{pe7}
\end{align}
In addition,
\begin{equation} \label{pe8} 
\| \Grad \vc{R}_\ep (\vc{h}_\ep (t) ) [ \Grad \Del^{-1} b(\vr_{\ep,f}) ] \cdot \vc{Y}_\ep	\|_{L^r(R^3)} \aleq |\vc{Y}_\ep|\ \| b(\vr_{\ep,f}) \|_{L^r(R^3)} 
\aleq |\vc{Y}_\ep|,\ 1 < r \leq \frac{\gamma}{\alpha}.
\end{equation}

Finally, 
\begin{align} 
\left\| \vc{R}_\ep (\vc{h}_\ep (t) ) \left[ \Grad \Del^{-1} \left[ \Big( b(\vr_{\ep, f}) - b'(\vr_{\ep, f}) \vr_{\ep, f} \Big) \Div \vue \right]  \right]	\right\|_{L^s(R^3)}
\aleq 1 \br \mbox{for all}\ \frac{3}{2} < s \leq \frac{3r}{3 - r},\ \frac{1}{2} + \frac{\alpha}{\gamma} = \frac{1}{r}. 
\label{pe9}
\end{align} 

Seeing that 
\begin{align}
(\vre \vue)_{\ep > 0} \ &\mbox{bounded in}\ L^\infty \left( 0,T; \left( L^{\frac{2 \gamma}{\gamma + 1}} + L^{\frac{2 \Gamma}{\Gamma + 1}} \right) (R^3; R^3) \right),\br 
(\vue)_{\ep > 0} \ &\mbox{bounded in}\ L^2(0,T; L^6(R^3; R^3)),\br  
(\vre)_{\ep > 0} \ &\mbox{bounded in}\ L^\infty \left( 0,T; \left( L^\gamma + L^\Gamma \right) (R^3) \right)
\nonumber
\end{align}
we may combine \eqref{pe7}, \eqref{pe9} to obtain 
\[
\int_0^T \psi \intRd{ \vre \vue \cdot \vc{R}_\ep (\vc{h}_\ep (t) ) [ \Grad \Del^{-1} \partial_t b(\vr_{\ep,f}) ]    } \dt
\ \mbox{bounded uniformly for}\ \ep \to 0.
\]

We conclude by estimating 
\begin{align}
\intRd{ \vre \vue \cdot \left[ \vc{R}_\ep (\vc{h}_\ep (t) ) [ \Grad \Grad \Del^{-1} b(\vr_{\ep,f})] \cdot \vc{Y}_\ep   - 
	\Grad \vc{R}_\ep (\vc{h}_\ep (t) ) [ \Grad \Del^{-1} b(\vr_{\ep,f}) ] \cdot \vc{Y}_\ep   \right]  }  \br 
= \int_{B_{2\ep}(\vc{Y}_\ep)}  \vre \vue \cdot \left[ \vc{R}_\ep (\vc{h}_\ep (t) ) [ \Grad \Grad \Del^{-1} b(\vr_{\ep,f})] \cdot \vc{Y}_\ep   - 
	\Grad \vc{R}_\ep (\vc{h}_\ep (t) ) [ \Grad \Del^{-1} b(\vr_{\ep,f}) ] \cdot \vc{Y}_\ep   \right]  \dx .
\label{pe10}
\end{align}
This integral can be decomposed as 
\begin{align}
\int_{B_{2\ep}(\vc{Y}_\ep)}  \vre \vue \cdot \left[ \vc{R}_\ep (\vc{h}_\ep (t) ) [ \Grad \Grad \Del^{-1} b(\vr_{\ep,f})] \cdot \vc{Y}_\ep   - 
	\Grad \vc{R}_\ep (\vc{h}_\ep (t) ) [ \Grad \Del^{-1} b(\vr_{\ep,f}) ] \cdot \vc{Y}_\ep   \right]  \dx  \br = 
\int_{B_{\ep}(\vc{Y}_\ep)}  \vr_{\ep,B} \vue \cdot \left[ \vc{R}_\ep (\vc{h}_\ep (t) ) [ \Grad \Grad \Del^{-1} b(\vr_{\ep,f})] \cdot \vc{Y}_\ep   - 
\Grad \vc{R}_\ep (\vc{h}_\ep (t) ) [ \Grad \Del^{-1} b(\vr_{\ep,f}) ] \cdot \vc{Y}_\ep   \right]  \dx  \br 
+ \int_{B_{2\ep}(\vc{Y}_\ep)}  \vr_{\ep,f} \vue \cdot \left[ \vc{R}_\ep (\vc{h}_\ep (t) ) [ \Grad \Grad \Del^{-1} b(\vr_{\ep,f})] \cdot \vc{Y}_\ep  - 
\Grad \vc{R}_\ep (\vc{h}_\ep (t) ) [ \Grad \Del^{-1} b(\vr_{\ep,f}) ] \cdot \vc{Y}_\ep   \right]  \dx.	
	\label{pe11}
\end{align}

Now, in accordance with \eqref{pe7}, \eqref{pe8},
\begin{align}
\left| \int_{B_{2\ep}(\vc{Y}_\ep)}  \vr_{\ep,f} \vue \cdot \left[ \vc{R}_\ep (\vc{h}_\ep (t) ) [ \Grad \Grad \Del^{-1} b(\vr_{\ep,f})] \cdot \vc{Y}_\ep   - 
\Grad \vc{R}_\ep (\vc{h}_\ep (t) ) [ \Grad \Del^{-1} b(\vr_{\ep,f}) ] \cdot \vc{Y}_\ep   \right]  \dx \right| \br 
\aleq \| \vr_{\ep, f} \|_{L^\gamma(R^3)} \| \vue \|_{L^6(R^3; R^3)} \| b(\vr_{\ep,f}) \|_{L^{\frac{\gamma}{\alpha}}(R^3)} |\vc{Y}_\ep | \ep^{3s},\ 
{s} = 1 - \frac{1}{\gamma} - \frac{\alpha}{\gamma} - \frac{1}{6} 
\label{pe12}
\end{align}
In view of \eqref{U2}, 
\[
|\vc{Y}_{\ep}| \aleq \ep^{\frac{1}{2} (\underline{\beta} - 3)} ,\ \underline{\beta} > 2 \frac{3 - \gamma}{\gamma}. 
\]
Consequently, 
\[
|\vc{Y}_{\ep}| \ep^{3s} \aleq \ep^{ \frac{1}{2} \underline{\beta} - \frac{3}{2} + 3 - \frac{3}{\gamma} - \frac{3\alpha}{\gamma} - \frac{1}{2} },\ \mbox{with}\ 
\frac{1}{2} \underline{\beta} - \frac{3}{2} + 3 - \frac{3}{\gamma} - \frac{3\alpha}{\gamma} - \frac{1}{2} = \frac{1}{2}\underline{\beta} - \frac{3 - \gamma}{\gamma} 
 - \frac{3 \alpha}{\gamma} > 0
\]
as long as $\alpha > 0$ is small enough.

Finally,
\begin{align}
\left| \int_{B_{\ep}(\vc{Y}_\ep)}  \vr_{\ep,B} \vue \cdot \left[ \vc{R}_\ep (\vc{h}_\ep (t) ) [ \Grad \Grad \Del^{-1} b(\vr_{\ep,f})] \cdot \vc{Y}_\ep   - 
\Grad \vc{R}_\ep (\vc{h}_\ep (t) ) [ \Grad \Del^{-1} b(\vr_{\ep,f}) ] \cdot \vc{Y}_\ep   \right]  \dx \right| \br 
\aleq \vr_{\ep, B} |\vc{Y}_\ep | \| \vue \|_{L^6(R^3; R^3)} \| b(\vr_{\ep,f}) \|_{L^{\frac{\gamma}{\alpha}}(R^3)} \ep^{3 \left( 1 - \frac{1}{6} - \frac{\alpha}{\gamma} \right)}, 
\nonumber
\end{align}
where
\[
\vr_{\ep, B} |\vc{Y}_\ep | \aleq \ep^{-\frac{3}{2}} \ep^{- \frac{\Ov{\beta}}{2} },\ \Ov{\beta} < 2. 
\]
Thus if $\alpha > 0$ is small enough, we get 
\[
\vr_{\ep, B} |\vc{Y}_\ep | \ep^{3 \left( 1 - \frac{1}{6} - \frac{\alpha}{\gamma} \right)} \to 0 \ \mbox{as}\ \ep \to 0.
\]

We have shown that the integrals $I_{5,\ep}$ remain bounded as $\ep \to 0$, which completes the proof of the pressure estimates claimed in \eqref{equip}.

\section{Convergence}
\label{c}

Our ultimate goal is to perform the limit in the momentum equation \eqref{P13}. To this end, we consider a smooth function 
\[
\bfphi \in C^k_c([0,T) \times R^3; R^3),\ k \geq 2 
\ \mbox{and its restriction}\ E_\ep (\vc{h}_{\ep}(t)))[ \bfphi (t, \cdot)], 
\]
where the latter is an eligible test function in \eqref{P13}.

\subsection{Time derivative}

We start with the time derivative 
\[
\int_0^T \intRd{ \vre \vue (t, \cdot) \cdot \partial_t \left( E_\ep (\vc{h}_{\ep}(t)))[ \bfphi (t, \cdot)] \right) } \dt.
\]
By virtue of formula \eqref{T11}, 
\begin{align}
\int_0^T \intRd{ \vre \vue (t, \cdot) \cdot \partial_t \left( E_\ep (\vc{h}_{\ep}(t)))[ \bfphi (t, \cdot)] \right) } \dt = 
\int_0^T \intRd{ \vre \vue (t, \cdot) \cdot  \left( E_\ep (\vc{h}_{\ep}(t)))[ \partial_t \bfphi (t, \cdot)] \right) } \dt \br 
+ \int_0^T \intRd{ \vre \vue (t, \cdot) \cdot \Big[ E_\ep (\vc{h}_{\ep}(t))[\Grad \bfphi(t, \cdot) ] \cdot \vc{Y}_\ep (t)
	- \Grad E_\ep (\vc{h}_{\ep}(t)) [\bfphi (t, \cdot)] \cdot \vc{Y}_\ep (t) \Big] } \dt.
\label{c1}
	\end{align}
In accordance with the convergence \eqref{U8b} and the estimate \eqref{lpe}, we get:  
\begin{equation} \label{c2}
\int_0^T \intRd{ \vre \vue (t, \cdot) \cdot  \left( E_\ep (\vc{h}_{\ep}(t)))[ \partial_t \bfphi (t, \cdot)] \right) } \dt \to 
\int_0^T \intRd{ \vr \vu \cdot  \partial_t \bfphi (t, \cdot) } \dt \ \mbox{as}\ \ep \to 0.
\end{equation}

As for the remaining integral in \eqref{c1}, we use \eqref{TT7} obtaining
\begin{align}
&\intRd{ \vre \vue (t, \cdot) \cdot \Big[ E_\ep (\vc{h}_{\ep}(t))[\Grad \bfphi(t, \cdot) ] \cdot \vc{Y}_\ep (t)
	- \Grad E_\ep (\vc{h}_{\ep}(t)) [\bfphi (t, \cdot)] \cdot \vc{Y}_\ep (t) \Big] } \br & \quad = 
\int_{B_{\frac{7}{4} \ep}(\vc{h}_\ep) } \vre \vue (t, \cdot) \cdot \Big[ E_\ep (\vc{h}_{\ep}(t))[\Grad \bfphi(t, \cdot) ] \cdot \vc{Y}_\ep (t)
- \Grad E_\ep (\vc{h}_{\ep}(t)) [\bfphi (t, \cdot)] \cdot \vc{Y}_\ep (t) \Big] \dx \br &\quad =
\int_{B_{\ep}(\vc{h}_\ep)}\vr_{\ep, B} \vue (t, \cdot) \cdot \Big[ E_\ep (\vc{h}_{\ep}(t))[\Grad \bfphi(t, \cdot) ] \cdot \vc{Y}_\ep (t)
- \Grad E_\ep (\vc{h}_{\ep}(t)) [\bfphi (t, \cdot)] \cdot \vc{Y}_\ep (t) \Big] \dx \br &\quad + 
\int_{ B_{\frac{7}{4} \ep}(\vc{h}_\ep) } \vr_{\ep, f} \vue (t, \cdot) \cdot \Big[ E_\ep (\vc{h}_{\ep}(t))[\Grad \bfphi(t, \cdot) ] \cdot \vc{Y}_\ep (t)
- \Grad E_\ep (\vc{h}_{\ep}(t)) [\bfphi (t, \cdot)] \cdot \vc{Y}_\ep (t) \Big] \dx.
\label{c3}
\end{align}
By virtue of \eqref{TT7}, \eqref{U2}, and hypothesis \eqref{P10},
\begin{align} 
&\left| \int_{B_{\ep}(\vc{h}_\ep)} \vr_{\ep, B} \vue (t, \cdot) \cdot \Big[ E_\ep (\vc{h}_{\ep}(t))[\Grad \bfphi(t, \cdot) ] \cdot \vc{Y}_\ep (t)
- \Grad E_\ep (\vc{h}_{\ep}(t)) [\bfphi (t, \cdot)] \cdot \vc{Y}_\ep (t) \Big] \dx \right| \br &\aleq 
\quad \int_{B_{\ep}(\vc{h}_\ep)} \sqrt{ \vr_{\ep, B} } | \vue (t, \cdot) | \ep^{-\frac{3}{2}} \dx .
\label{c4}
\end{align}
Thus, in view of the uniform bounds \eqref{U1f}, 
we conclude
\begin{equation} \label{c6}
\int_0^T \int_{B_{\ep}(\vc{h}_\ep)} \vr_{\ep, B} \vue (t, \cdot) \cdot \Big[ E_\ep (\vc{h}_{\ep}(t))[\Grad \bfphi(t, \cdot) ] \cdot \vc{Y}_\ep (t)
- \Grad E_\ep (\vc{h}_{\ep}(t)) [\bfphi (t, \cdot)] \cdot \vc{Y}_\ep (t) \Big] \dx \dt \to 0  
\end{equation}
as $\ep \to 0$.

The second integral on the right--hand side of \eqref{c3} can be handled by using \eqref{lpe} and \eqref{T7}:
\begin{align}
&\left\| \vr_{\ep, f} \vue (t, \cdot) \cdot \Big[ E_\ep (\vc{h}_{\ep}(t))[\Grad \bfphi(t, \cdot) ] \cdot \vc{Y}_\ep (t)
- \Grad E_\ep (\vc{h}_{\ep}(t)) [\bfphi (t, \cdot)] \cdot \vc{Y}_\ep (t) \Big] \right\|_{L^1(R^3)} \br &\quad \aleq
\| \vr_{\ep, f}(t ,\cdot) \|_{L^\gamma (R^3)} \| \vue (t, \cdot) \|_{L^6(R^3)} \| \vc{Y}_\ep \mathds{1}_{B_{2\ep}(\vc{h}_\ep)} \|_{L^s(R^3)}\|\Grad \bfphi\|_{L^{\infty}(R^3)},\ 
\frac{1}{\gamma} + \frac{1}{6} + \frac{1}{s} = 1, 	
	\label{c7}
	\end{align}
Moreover, by virtue of \eqref{U2},  
\begin{equation} \label{c8}
|\vc{Y}_\ep (t) | \leq \frac{1}{ \sqrt{\vr_{\ep, B}}} \ep^{\frac{-3}{2}} \aleq \ep^{\frac{\underline{\beta} - 3}{2}}.
\end{equation}
Thus, it follows from hypothesis \eqref{P10} and a direct manipulation \begin{align}
&\int_0^T \left| \int_{ B_{\frac{7}{4} \ep}(\vc{h}_\ep) }\vr_{\ep, f} \vue (t, \cdot) \cdot \Big[ E_\ep (\vc{h}_{\ep}(t))[\Grad \bfphi(t, \cdot) ] \cdot \vc{Y}_\ep (t)
- \Grad E_\ep (\vc{h}_{\ep}(t)) [\bfphi (t, \cdot)] \cdot \vc{Y}_\ep (t) \Big] \dx \dt \right| \br &\quad \to 0 
\ \mbox{as}\ \ep \to 0.
\label{c9}
\end{align}

Summing up \eqref{c2}, \eqref{c6}, and \eqref{c9} we conclude 
\begin{equation} \label{c10}
\int_0^T \intRd{ \vre \vue (t, \cdot) \cdot \partial_t \left( E_\ep (\vc{h}_{\ep}(t)))[ \bfphi (t, \cdot)] \right) } \dt \to 
\int_0^T  \intRd{ \vr \vu \cdot  \partial_t \bfphi (t, \cdot) } \dt \ \mbox{as}\ \ep \to 0.
	\end{equation}

\subsection{Convective term and the viscous stress}

Repeating the arguments of the previous section, we easily establish
\begin{align} 
\int_0^T & \intRd{ \vre \vue \otimes \vue (t, \cdot) : \Grad \left( E_\ep (\vc{h}_{\ep}(t)))[ \bfphi (t, \cdot)] \right) } \dt \br 
&=  \int_0^T \intRd{ \vr_{\ep,f} \vue \otimes \vue (t, \cdot) : \Grad \left( E_\ep (\vc{h}_{\ep}(t)))[ \bfphi (t, \cdot)] \right) } \dt 
\br &\to \int_0^T \intRd{ \Ov{\vr \vu \otimes \vu } : \Grad \bfphi } \dt \ \mbox{as}\ \ep \to 0,
\label{c11}	
	\end{align}
and 
\begin{align} 
	\int_0^T & \intRd{ \mathbb{S} (\Grad \vue) : \Grad \left( E_\ep (\vc{h}_{\ep}(t)))[ \bfphi (t, \cdot)] \right) } \dt \br 
	&\to \int_0^T \intRd{ \mathbb{S} (\Grad \vu) : \Grad \bfphi } \dt \ \mbox{as}\ \ep \to 0.
	\label{c12}	
\end{align}
Here, in view of the uniform bounds \eqref{inter}, 
\[
\vr_{\ep, f} \vue \otimes \vue \to 
\Ov{\vr \vu \otimes \vu } \ \mbox{weakly in}\ L^p_{\rm loc}([0,T] \times R^3; R^{3 \times 3}) \ \mbox{for some}\ p > 1.
\]

Our ultimate goal in the section is to prove the identity
\begin{equation} \label{c13}
\Ov{\vr \vu \otimes \vu } = \vr \vu \otimes \vu.
\end{equation}	
To this end, consider 
\[
\bfphi = E_\ep (\vc{h}_{\ep}(t))) [\psi(t) \phi(\cdot) ] ,\ \psi \in C^1_c(0,T),\ \phi \in C^1_c(R^3; R^3)
\]
as a test function in the momentum equation \eqref{P13}. We easily compute
\begin{align}
	\int_0^T &\intRd{ \Big[ \vre \vue \cdot \partial_t \left( E_\ep (\vc{h}_{\ep}(t))) [\psi(t) \phi(\cdot) ] \right) } \dt = - \int_0^T \psi (t) \intRd{ \vre \vue \otimes \vue : \Grad 
E_\ep (\vc{h}_{\ep}(t))) [\phi(\cdot) ] } \dt \br		
		&- \int_0^T \psi(t) \intRd{  p(\vre) \Div E_\ep (\vc{h}_{\ep}(t))) [\phi(\cdot) ]  } \dt
\br &+ \int_0^T \psi(t) \intRd{ \mathbb{S}(\Grad \vue) : \Grad E_\ep (\vc{h}_{\ep}(t))) [ \phi(\cdot) ] } \dt. 
\label{c14}
\end{align}
Moreover, by virtue of formula \eqref{T11},
\begin{align}
\int_0^T &\intRd{ \Big[ \vre \vue \cdot \partial_t \left( E_\ep (\vc{h}_{\ep}(t))) [\psi(t) \phi(\cdot) ] \right) } \dt 
\br &= \int_0^T \partial_t \psi(t) \intO{ \vre \vue \cdot E_\ep (\vc{h}_{\ep}(t)))[ \phi (\cdot) ] } \dt \br & + 
\int_0^T \psi (t) \intRd{ \vre \vue \cdot E_\ep (\vc{h}_{\ep}(t))[\Grad \phi(\cdot) ] \cdot \vc{Y}_\ep (t) } \dt \br 
&- \int_0^T \psi(t) \intRd{ \vre \vue \cdot \Grad E_\ep (\vc{h}_{\ep}(t)) [\phi (\cdot)] \cdot \vc{Y}_\ep (t) } \dt.
\label{c15}
\end{align}
If $\phi$ is smooth ($C^1_c(R^3; R^3)$), the last two integrals in \eqref{c15} can be handled exactly as their counterpart in 
\eqref{c6}, \eqref{c9}, specifically,
\begin{align} 
&\left\| \intRd{ \vre \vue \cdot \left( E_\ep (\vc{h}_{\ep}(t))[\Grad \phi(\cdot) ] -   \Grad E_\ep (\vc{h}_{\ep}(t)) [\phi (\cdot)] \right) \cdot\vc{Y}_\ep (t) }\right\|_{L^2(0,T)} \br 
&\quad \leq
c(\| \phi \|_{C^1}) \ \mbox{independently of}\ \ep \to 0. 	
\label{c16}
\end{align}
Combining \eqref{c14}--\eqref{c16} we may infer that the function
\begin{equation} \label{c17}
t \in [0,T] \mapsto \intRd{ \vre \vue (t,\cdot) \cdot E_\ep (\vc{h}_{\ep}(t))) [\phi ] } ,\ \phi \in C^1_c(R^3; R^3),
\end{equation}
is H\"older continuous with a positive exponent and norm depending solely on $\| \phi \|_{C^1(\Ov{\Omega}; R^d)}$.

Finally, by virtue of the error estimates \eqref{TT7}, 
\[
E_\ep (\vc{h}_{\ep}(t))) [\phi ] \to \phi \ \mbox{in}\ W^{1,2}(R^3; R^3) \ \mbox{for any}\ \phi \in W^{1,2}(R^3; R^3)\ \mbox{uniformly in}\ t \in (0,T),
\]
and we deduce from \eqref{c17} that 
\begin{equation} \label{c18}
	\vre \vue \ \mbox{precompact in}\ L^2(0,T; W^{-1,2}_{\rm loc}(R^3; R^3)), 
	\end{equation}
which, together with \eqref{U3} yields \eqref{c13}.

\subsection{The pressure and strong convergence of the density}

In view of the pressure estimates \eqref{equip}, it is easy to establish the limit 
\[
\int_0^T \intRd{ p(\vr_{\ep,f}) \Div E_\ep (\vc{h}_\ep) [\bfphi ] } \dt \to \int_0^T \intRd{ \Ov{p(\vr)} \Div \bfphi  } \dt,
\]
where $\Ov{p(\vr)}$ stands for a weak limit of the sequence $(p(\vr_{\ep,f}))_{\ep > 0}$. 

Thus the remaining issue is to establish the equality 
\[
p(\vr) = \Ov{p(\vr)} 
\]
which is the standard and nowadays well understood problem in the theory of compressible fluids, see e.g. \cite{EF70}, Lions \cite{LI4}.
The proof requires $\bfphi = \psi(t) \phi(x) \Grad \Del^{-1}[ b(\vr_{\ep,f})]$ to be used as test functions in the momentum balance, where $b$ is a bounded functions and 
$\psi \in C^1_c(0,T)$, $\phi \in C^1_c(R^3)$.
In the present setting, similarly to the above, we use the quantity 
\[
\bfphi = \psi E_\ep (\vc{h}_\ep) \left[ \phi \Grad \Del^{-1}[ b(\vr_{\ep,f})] \right],
\]
which is a legal test function for the momentum balance \eqref{P13}. As shown above, the resulting error terms vanish in the asymptotic limit $\ep \to 0$ and the proof of 
the strong convergence of the density is therefore the same as in the fluid without moving objects. Thus exactly the same method as in \cite[Chapter 6]{EF70} can be used to complete the proof of Theorem 
\ref{PT1}.



\def\cprime{$'$} \def\ocirc#1{\ifmmode\setbox0=\hbox{$#1$}\dimen0=\ht0
	\advance\dimen0 by1pt\rlap{\hbox to\wd0{\hss\raise\dimen0
			\hbox{\hskip.2em$\scriptscriptstyle\circ$}\hss}}#1\else {\accent"17 #1}\fi}

\end{document}